\documentclass[a4papert]{article}
\usepackage{amsfonts}
\usepackage{amssymb}
\usepackage{mathrsfs}
\usepackage{float}
\usepackage{amsmath}
\usepackage{amsfonts,amssymb}
\usepackage{graphicx,color}
\usepackage{mathrsfs}

\newtheorem{Theorem}{Theorem}[section]

\newtheorem{Lemma}{Lemma}[section]

\newtheorem{Remark}{Remark}[section]

\def\2{{I \hskip -1.0mm I}}
\def\3{{I \hskip -1.0mm I\hskip -1.0mm I}}
\def\4{{I \hskip -0.9mm V}}
\def\6{{V \hskip -1.35mm I}}

\date{ }

\title{Generalized Ricci flow I: Local existence and
uniqueness}

\author{Chun-Lei He\footnote{Department of Mathematics, Shanghai Jiao Tong University,
Shanghai 200240, China; }\,\,, $\quad$ Sen Hu\footnote{Department
of Mathematics, University of Science and Technology of China,
Hefei 230026, China;}\,\,, $\quad$De-Xing Kong\footnote{Center of
Mathematical Sciences, Zhejiang University, Hangzhou 310027,
China;} $\quad$ and $\quad$ Kefeng Liu\footnote{ Department of
Mathematics, University of California at Los Angeles, CA 90095,
USA.}\\}
\begin{document}
\maketitle
\begin{abstract}
In this paper we investigate a kind of generalized Ricci flow
which possesses a gradient form. We study the monotonicity of the
given function under the generalized Ricci flow and prove that the
related system of partial differential equations are strictly and
uniformly parabolic. Based on this, we show that the generalized
Ricci flow defined on a $n$-dimensional compact Riemannian
manifold admits a unique short-time smooth solution. Moreover, we
also derive the evolution equations for the curvatures, which play
an important role in our future study.
\end{abstract}


\section{Introduction}
 In the early eighties R. Hamilton introduced the Ricci
flow to construct canonical metrics for some manifolds. Since then
many mathematicians, including Hamilton, Yau, Perelman and others,
developed many tools and techniques to study the Ricci flow. The
latest developments confirmed that the Ricci flow approach is very
powerful in the study of three-manifolds. In fact, a complete
proof of Poincare's conjecture and Thurston's geometrization
conjecture has been offered in Cao-Zhu's paper \cite{Cao} and
others after Perelman's breakthrough.

It is useful to observe that, in Perelman's work \cite{Perelman},
a key step is to introduce a functional for a metric $g$ and a
function $f$ on a manifold $M$
$$ W(g, f) = \int_{M^{3}}d^{3}x \sqrt{g} e^{-f} (R + |\nabla
f|^{2}).$$ The variation of this functional generates a gradient
flow which is a system of partial differential equations
$$\dot{g}_{ij} = -2 (R_{ij} + \nabla_{i} \nabla_{j} f),$$
$$ \dot{f} = - (R + \triangle f).$$
If we fix a measure for the conformal class of metrics $e^{f}
ds^{2}$ of a metric, i.e., let $dm = e^{-f} dV$ be fixed, then we
get back to the original Ricci flow after we apply a
transformation of diffeomorphism generated by the vector field
$\nabla_{i} f$ to the metric. In this way, we express the Ricci
flow as a gradient flow. Dynamics of a gradient flow is much
easier to handle. The functional generating the flow gives a
monotone functional along the orbit of the flow automatically. If
the flow exists for all time, then it shall flow to a critical
point which leads to the existence of a canonical metric. Even for
a flow which does not exist for all time, the generating
functional helps very much in the analysis of singularities.

Perelman's above idea came from physics. Ricci flow arises as the
first order approximation of the renormalization flow of a sigma
model. Since there are many kinds of sigma models, it would be
interesting to try some other models. Indeed such a generalization
was made by physicists in \cite{Thurston}. For a three-manifold
$M^{3}$, they proposed to add a $U(1)$ gauge field with potential
1-form $A$ and field strength $F$ which are coupled as a {\it
Maxwell-Chern-Simons theory}. The corresponding action given by
\cite{Gegenberg1} or \cite{Gegenberg2} reads
$$S = \int_{M} d^{3}x \sqrt{g} e^{- f} ( -\chi + R + |\nabla
f|^{2}) - \frac{1}{2}e^{- f} H \wedge * H -e^{- f} F \wedge * F
.$$ The $U(1)$ gauge field $A$ is a one-form potential whose field
strength $F=dA.$ The Wess-Zumino field $B$ is a two-form potential
whose field strength $H = dB ,$  $f$ is a dilaton. In their paper,
they find that Thurston's eight geometries appear as critical
points of the above functional. Furthermore they show that there
are no other critical points. So basically critical points of the
above functional are eight geometries of Thurston. They also
propose to study the gradient flow of the functional $S$ as a
generalization of the Ricci flow. Unfortunately, they modify the
gradient flow in a way to change sign for the variable of gauge
fields. Although the modified flow shares the same set of critical
points they lost the important monotone property (along an orbit).

 In addition, we are also able to consider a flow for a
similar functional for a four-dimension manifold
$$S_1 = \int_{M} d^{4}x \sqrt{g} e^{- f} ( \chi + R + |\nabla f
|^{2}) - \frac{1}{2}e^{- f} H \wedge * H - e^{- f} F \wedge
* F + \frac{e}{2} F \wedge F$$where $e$ is the Euler number $e(\eta)$
 of the bundle $\eta.$ The corresponding flow is given by
$$\begin{cases} \dfrac{\partial g_{ij}}{\partial
t}=-2[R_{ij}+\nabla_i\nabla_j f-\dfrac{1}{4}H_{ikl}H_j^{~kl}-F_i^{~k}F_{jk}],\vspace{2mm}\\
\dfrac{\partial B_{ij}}{\partial
t}=e^{f}\nabla_k(e^{-f}H^k_{~ij}),\vspace{2mm}\\
\dfrac{\partial A_i}{\partial
t}=-e^{f}\nabla_k(e^{-f}F^{~k}_{i}),\vspace{2mm}\\
\dfrac{\partial f}{\partial t}=\chi-2R-3\triangle f+|\nabla
f|^2+\dfrac{1}{3}H^2+\dfrac{3}{2}F^2.
\end{cases}$$

The generalization to four-manifolds is probably more interesting.
It may offer a systematic way to study four-manifolds.

The success of studying three-manifolds relies on a program
proposed by Thurston, i.e., his geometrization conjecture. He
conjectures and proves for several large classes of
three-manifolds, that every three-manifold can be decomposed into
pieces of three-manifolds of canonical metrics, i.e., those
manifolds carrying one of the eight geometries of Thurston.

For four-dimension manifolds the critical points of $S_1$ might
play a similar role as building blocks of smooth four-dimension
manifolds. It would be interesting to study those critical points
and to study what other four manifolds one can get by performing
surgeries and gluing on those manifolds. We shall address this
problem in the future.

 As a first
step, we shall show that the flow does exist. We shall also prove
that the modified system of partial differential equations are
strictly and uniformly parabolic.

The paper is organized as follows. Section 2 is devoted to the
proof of local existences and uniqueness. In section 3 we study
the monotonicity of $S$ under the modified flow. In Section 4 we
investigate the equations for the critical points of $S$ and point
out that fields $F$ and $H$ do not provide any help for the case
of compact manifold but maybe play an important role for the
noncompact case. In Section 5, we derive the evolution equations
for the curvatures, which play an important role in our future
study.

\section{Local Existences and Uniqueness}

In this section, we mainly establish the short-time existence and
uniqueness result for the gradient flow (\ref{14}), (\ref{15}) and
(\ref{16}) on a compact 3-dimensional manifold $M$. It is known
that the gradient flow (\ref{14}), (\ref{15}) and (\ref{16})  is a
system of second order nonlinear weakly parabolic partial
differential equations. By the proof of the local existence and
uniqueness of the Ricci flow (for example see \cite{Cao} \cite{D},
), we can obtain a modified evolution equations by the
diffeomorphism $\varphi$ of $M$, which is a strictly parabolic
system. Then, by the standard theory of parabolic equations, the
modified evolution equations has a uniqueness solution.

Let us choose a normal coordinate $\{x^i\}$ around a fixed point
$x\in M$ such that $\dfrac{\partial g_{ij}}{\partial x^k}=0$ and
$g_{ij}(p)=\delta_{ij}.$

\begin{Theorem} {\bf (Local existences and uniqueness)} Let $(M,
g_{ij}(x)) $ be a three-dimensional compact Riemannian manifold.
Then there exists a constant $T
> 0$ such that the evolution equations
\begin{equation}\label{21}\begin{cases} \dfrac{\partial g_{ij}}{\partial
t}=-2[R_{ij}-\dfrac14H_{ikl}H_j^{~kl}-F_i^{~k}F_{jk}],\vspace{2mm}\\
\dfrac{\partial A_i}{\partial t}=-\nabla_k {F}_{i}^{~k},\\
\dfrac{\partial B_{ij}}{\partial
t}=\nabla_k{H}^{k}_{~ij}.\vspace{2mm}\end{cases}\end{equation}
 has a unique smooth solution on $M\times[0,T)$ for
every initial fields.
\end{Theorem}

\begin{Lemma}  For each gauge equivalent class of a gauge field
$A$, there exists an $A^{'}$ such that $d (* A^{'}) =
0$.\end{Lemma}

The lemma can be proved by the Hodge decomposition.

\noindent{\bf Proof.} For each one-form $A$, by the Hodge
decomposition, there exists an one-form $A_{0}$, a function
$\alpha$ and a two-form $\beta$ such that
$$A = A_{0} + d \alpha + d^{*} \beta,$$
$$d A_{0} = 0, ~~~d^{*}A_{0} = 0.$$
Let $A'=A-d\alpha$. $A'$ is in the same gauge equivalent class of
$A$. Since $d(*A_0)=0, \,d(*d^*\beta)=0,$ then we have $d(*A')=0$.
$\qquad\Box$

\begin{Lemma}  The differential operator of the right hand of (\ref{15})
with respect to the gauge equivalent class of a gauge field $A$ is
uniformly elliptic. \end{Lemma}

\noindent{\bf Proof.} Let $A=A_idx^i$ be a gauge field. By Lemma
4.1 we can choose an $A^{'}$ in the gauge equivalent class of $A$
such that $d(* A^{'}) = 0$. We still denote $A'$ as $A$ . Since
$d(*A)=0,$ we have $dd^*A=0 ,$ then
$\sum\limits_{k=1}^{3}\dfrac{\partial^2 A_k}{\partial x^k\partial
x^i}=0, \forall~ i=1,\cdots,3.$ Noting that $F = dA$ and
$F_{ij}=\dfrac{\partial A_j}{\partial x^i}-\dfrac{\partial
A_i}{\partial x^j},$ We have
\begin{eqnarray*}\dfrac{\partial A_i}{\partial t }&
=&-\nabla_k {F}_{i}^{~k} =- \nabla_{k}(g^{kl}F_{il})=
-g^{kl}(\dfrac{\partial^2A_l}{\partial x^k\partial
x^i}-\dfrac{\partial^2A_i}{\partial x^k\partial x^l})
=g^{kl}\dfrac{\partial^2A_i}{\partial x^k\partial
x^l}.\end{eqnarray*} The right hand side of above equation is
clearly elliptic at point $x$. If we apply a diffeomorphism to the
metric it won't change the positivity property of the second order
operator of the right hand side.$\qquad\Box$

Now let us consider the equation for $B_{ij}$.

\begin{Lemma}  For each gauge equivalent class of a B-field B,
i.e., a two-form B on $M$, there exists a $B^{'}$ such that $d( *
B^{'}) = 0$.\end{Lemma}

\noindent{\bf Proof.} Again we use the Hodge decomposition. For a
two-form $B$, there exist a one-form $\alpha$, a two-form $B_{0}$
and a three-form $\beta$ such that
$$B = B_{0} + d \alpha + d^{*} \beta,$$
$$d B_{0} = 0, d^{*} B_{0} = 0.$$

Let $B^{'} = B - d \alpha.$ Since $ B^{'}$ is in the same gauge
equivalent class of $B$,  we have $d(*B^{'}) = 0.$$\qquad\Box$

\begin{Lemma}  The differential operator of the right hand side of
(\ref{16}) with respect to the gauge equivalent class of a
$B$-field $B$ is uniformly elliptic.\end{Lemma}

\noindent{\bf Proof.} Let us consider the equation for $B$-field.
Without loss of generality, we assume $d(* B) = 0.$ Thus $dd^* B=
0.$ Then $\sum\limits_{k=1}^{3}(\dfrac{\partial^2B_{ki}}{\partial
x^k\partial x^j}+\dfrac{\partial^2B_{jk}}{\partial x^k\partial
x^i})=0, \forall~i,j=1,\cdots,3.$ We have
\begin{eqnarray*}
\dfrac{\partial B_{ij}}{\partial t} =
\nabla_kH^k_{~ij}=g^{kl}(\dfrac{\partial^2 B_{ij}}{\partial
x^k\partial x^l}+ \dfrac{\partial^2 B_{jl}}{\partial x^k\partial
x^i}+\dfrac{\partial^2 B_{li}}{\partial x^k\partial x^j})
=g^{kl}\dfrac{\partial^2 B_{ij}}{\partial x^k\partial x^l}.
\end{eqnarray*} The right hand side is clearly elliptic at the
point $x$. If we apply a diffeomorphism to the metric it does not
change the positivity property of the second order operator of the
right hand side.$\qquad\Box$

Suppose $\hat{g}_{ij}(x,t) $ is a solution of the equations
(\ref{21}), and $\varphi_t: M\rightarrow M$ is a family of
diffeomorphisms of $M$. Let
$$g_{ij}(x,t)=\varphi^*_t\hat{g}_{ij}(x,t),$$ where $\varphi^*_t$ is the pull-back operator of $\varphi_t.$ We
now want to find the evolution equations for the metric
$g_{ij}(x,t).$

Denote
$$y(x,t)=\varphi_t(x)=\{y^1(x,t), y^2(x,t),\cdots, y^n(x,t)\}$$ in
local coordinates. Then
$$g_{ij}(x,t)=\dfrac{\partial y^\alpha}{\partial x^i}\dfrac{\partial y^\beta}{\partial
x^j}\hat{g}_{\alpha\beta}(y,t)$$ and
\begin{eqnarray}\begin{array}{lll} \dfrac{\partial}{\partial
t}{g_{ij}(x,t)} &=& \dfrac{\partial}{\partial t}
\left[{\hat{g}_{\alpha\beta}(y,t)}\cdot\dfrac{\partial
y^\alpha}{\partial x^{i}}\cdot\dfrac{\partial y^\beta}{\partial
x^{j}}\right]\vspace{2mm} \\
&=& \dfrac{\partial y^\alpha}{\partial x^{i}}\dfrac{\partial
y^\beta}{\partial x^{j}}\dfrac{\partial}{\partial
t}\hat{g}_{\alpha\beta}(y,t)+\dfrac{\partial y^\alpha}{\partial
x^{i}}\dfrac{\partial y^\beta}{\partial x^{j}}\dfrac{\partial
y^\gamma}{\partial t}\dfrac{\partial}{\partial
y^\gamma}\hat{g}_{\alpha\beta}(y,t)\vspace{2mm}\\
&&+\hat{g}_{\alpha\beta}(y,t)\dfrac{\partial}{\partial
x^i}(\dfrac{\partial y^\alpha}{\partial t})\dfrac{\partial
y^\beta}{\partial x^j}+\hat{g}_{\alpha\beta}(y,t)\dfrac{\partial
y^\alpha}{\partial x^i}\dfrac{\partial}{\partial
x^j}(\dfrac{\partial y^\beta}{\partial
t})\,\,.\nonumber\end{array}\end{eqnarray} Since
\begin{eqnarray*}\begin{aligned}\dfrac{\partial y^\alpha}{\partial x^i}\dfrac{\partial y^\beta}
{\partial x^j}\dfrac{\partial y^\gamma}{\partial t}\dfrac{\partial
}{\partial y^\gamma}{\hat g_{\alpha\beta}}&=\dfrac{\partial
y^\alpha}{\partial x^i}\dfrac{\partial y^\beta} {\partial
x^j}\dfrac{\partial y^\gamma}{\partial
t}g_{kl}\dfrac{\partial}{\partial y^\gamma}(\dfrac{\partial
x^k}{\partial y^\alpha}\dfrac{\partial x^l}{\partial
y^\beta})\vspace{2mm}\\
&=\dfrac{\partial y^\beta}{\partial t}\dfrac{\partial^2
x^k}{\partial y^\alpha\partial y^\beta}\dfrac{\partial
y^\alpha}{\partial x^i}g_{jk}+\dfrac{\partial y^\alpha}{\partial
t}\dfrac{\partial^2 x^k}{\partial y^\alpha\partial
y^\beta}\dfrac{\partial y^\beta}{\partial
x^j}g_{ik}\,\,,\\
\Gamma^k_{jl}&=\dfrac{\partial y^\alpha}{\partial
x^j}\dfrac{\partial y^\beta}{\partial x^l}\dfrac {\partial
x^k}{\partial
y^\gamma}\hat{\Gamma}^\gamma_{\alpha\beta}+\dfrac{\partial
x^k}{\partial y^\alpha}\dfrac{\partial^2 y^\alpha}{\partial
x^j\partial x^l}\,\,,\\
\end{aligned}\end{eqnarray*}
 then
$$\dfrac{\partial y^\alpha}{\partial x^i}\dfrac{\partial y^\beta}
{\partial x^j}\hat H_{\alpha\rho\delta} \hat
H_\beta^{~\rho\delta}=H_{ikl}H_j^{~kl},~~~~\dfrac{\partial
y^\alpha}{\partial x^i}\dfrac{\partial y^\beta} {\partial x^j}\hat
F_\alpha^{~\rho} \hat F_{\beta\rho}=F_i^{~k}F_{jk}\,\,.$$

Therefore, in the normal coordinate, we have

\begin{eqnarray*}\begin{array}{lll}
\dfrac{\partial }{\partial
t}g_{ij}(x,t)&=&\dfrac{\partial}{\partial x^i}(\dfrac{\partial
y^\alpha}{\partial t})\dfrac{\partial y^\beta}{\partial
x^j}g_{kl}\dfrac{\partial x^k}{\partial y^\alpha}\dfrac{\partial
x^l}{\partial y^\beta}+\dfrac{\partial y^\alpha}{\partial
x^i}\dfrac{\partial}{\partial x^j}(\dfrac{\partial
y^\beta}{\partial t})g_{kl}\dfrac{\partial x^k}{\partial
y^\alpha}\dfrac{\partial x^l}{\partial y^\beta}+\dfrac{\partial
y^\alpha}{\partial t}\dfrac{\partial}{\partial
x^i}(\dfrac{\partial x^k}{\partial
y^\alpha})g_{jk}\vspace{2mm}\\
&&+\dfrac{\partial y^\beta}{\partial t}\dfrac{\partial }{\partial
x^j}(\dfrac{\partial x^k}{\partial y^\beta})g_{ik}+\dfrac{\partial
y^\alpha}{\partial x^i}\dfrac{\partial y^\beta} {\partial
x^j}\left[-2(\hat{R}_{\alpha\beta}-\dfrac{1}{4}\hat
H_{\alpha\rho\delta}
\hat H_{\beta}^{~\rho\delta}-\hat F_\alpha^{~\rho}\hat F_{\beta\rho})\right]\vspace{2mm}\\
&=&\dfrac{\partial}{\partial x^i}(\dfrac{\partial
y^\alpha}{\partial t}\dfrac{\partial x^k}{\partial
y^\alpha})g_{jk}+\dfrac{\partial}{\partial x^j}(\dfrac{\partial
y^\beta}{\partial t}\dfrac{\partial x^k}{\partial
y^\beta})g_{ik}-2R_{ij}+\dfrac{1}{2}H_{ikl}H_j^{~kl}+2F_i^{~k}F_{jk}
\vspace{2mm}\\
&=&-2R_{ij}+\nabla_i(\dfrac{\partial y^\alpha}{\partial
t}\dfrac{\partial x^k}{\partial
y^\alpha}g_{jk})+\nabla_j(\dfrac{\partial y^\beta}{\partial
t}\dfrac{\partial x^k}{\partial
y^\beta}g_{ik})\vspace{2mm}+\dfrac{1}{2}H_{ikl}H_j^{~kl}+2F_i^{~k}F_{jk}\,\,.\\
\end{array}\end{eqnarray*}
If we define $y(x,t)=\varphi_t(x)$ by the equations
\begin{equation}\label{17}\begin{cases}\dfrac{\partial y^\alpha}{\partial
t}=\dfrac{\partial y^\alpha}{\partial x^k}
(g^{jl}(\Gamma^k_{jl}-\tilde{\Gamma}^k_{jl})),\\
y^\alpha(x,0)=x^\alpha\end{cases}\end{equation}and
$V_i=g_{ik}g^{jl}(\Gamma^k_{jl}-\tilde{\Gamma}^k_{jl}),$ we get
the following evolution equations for the pull-back metric
\begin{equation}\label{18}\begin{cases}\dfrac{\partial }{\partial
t}g_{ij}(x,t)=-2R_{ij}+\nabla_iV_j+\nabla_jV_i+ \dfrac{1}{2}H_{ikl}H_j^{~kl}+2F_i^{~k}F_{jk},\\
g_{ij}(x,0)=\tilde{g}_{ij}(x),\\\end{cases}\end{equation}where $
\tilde g_{ij}(x)$
 is the initial metric
and $\tilde{\Gamma}^k_{jl}$ is the connection of the initial
metric. The initial value problem (\ref{17}) can be rewritten as
\begin{equation}\label{19}\begin{cases}\dfrac{\partial y^\alpha}{\partial
t}=g^{jl}(\dfrac{\partial^2 y^\alpha}{\partial x^j\partial x^l}
+\dfrac{\partial y^\beta}{\partial x^j}\dfrac{\partial
y^\gamma}{\partial
x^l}\hat{\Gamma}^\alpha_{\beta\gamma}-\dfrac{\partial
y^\alpha}{\partial x^k}\tilde {\Gamma}^k_{jl}
)\,\,,\\
y^\alpha(x,0)=x^\alpha.\end{cases}\end{equation} Equation
(\ref{19}) is clearly a strictly parabolic system. Then, we have
\begin{eqnarray*}\begin{array}{lll}\dfrac{\partial }{\partial t}g_{ij}(x,t)&=&\dfrac{\partial}{\partial x^i}
\left\{g^{kl}\dfrac{\partial g_{kl}}{\partial
x^j}\right\}-\dfrac{\partial}{\partial
x^k}\left\{g^{kl}\Large(\dfrac{\partial g_{jl}}{\partial
x^i}+\dfrac{\partial g_{il}}{\partial x^j}-\dfrac{\partial
g_{ij}}{\partial
x^l}\Large)\right\}\vspace{2mm}\\
&&+\dfrac{\partial}{\partial
x^i}\left\{g_{jk}g^{pq}\frac12g^{km}(\dfrac{\partial
g_{mq}}{\partial x^p}+\dfrac{\partial g_{mp}}{\partial
x^q}-\dfrac{\partial
g_{pq}}{\partial x^m})\right\}\vspace{2mm}\\
&&+\dfrac{\partial}{\partial
x^j}\left\{(g_{ik}g^{pq}\frac12g^{km}(\dfrac{\partial
g_{mq}}{\partial x^p}+\dfrac{\partial g_{mp}}{\partial
x^q}-\dfrac{\partial g_{pq}}{\partial x^m})\right\}+
\dfrac{1}{2}H_{ikl}H_j^{~kl}+2F_i^{~k}F_{jk}\vspace{2mm}\\
&=&g^{kl}\dfrac{\partial^2g_{ij}}{\partial x^k\partial x^l}+
\dfrac{1}{2}H_{ikl}H_j^{~kl}+2F_i^{~k}F_{jk}\,\,.\\\end{array}\end{eqnarray*}
 As a result, from the original equations, we can obtain
\begin{equation}\label{20}\begin{cases}\dfrac{\partial g_{ij}(x,t)}{\partial
t}=g^{kl}\dfrac {\partial^2g_{ij}}{\partial x^k\partial x^l}+
\dfrac{1}{2}H_{ikl}H_j^{~kl}+2F_i^{~k}F_{jk},\vspace{2mm}\\
\dfrac{\partial A_i}{\partial t
}=g^{kl}\dfrac{\partial^2A_i}{\partial x^k\partial
x^l}\vspace{2mm},\\
\dfrac{\partial B_{ij}}{\partial t }
=g^{kl}\dfrac{\partial^2 B_{ij}}{\partial x^k\partial x^l}.\\
\end{cases}\end{equation}
Let
\begin{equation*}\begin{aligned}&u_1=g_{11},u_2=g_{12},u_3=g_{13},u_4=g_{22},u_5=g_{23},u_6=g_{33},\\
&u_7=A_{1}, u_8=A_{2},
u_9=A_{3},u_{10}=B_{12},u_{11}=B_{13},u_{12}=B_{23}\,\,.\end{aligned}\end{equation*}
The above equations can be rewritten as the following
form$$\dfrac{\partial u_i}{\partial
t}=\sum_{j\,k\,l}a_{ikjl}\dfrac{\partial^2u_j}{\partial
x^k\partial x^l}+(\text lower \;order\; terms)\quad (k,l=1,2,3;\;
i,j=1,2,\cdots,12),$$ in which
$$a_{ikjl}=g^{kl}\;(j=i),~~~~ a_{ikjl}=0\;(j\neq i)~~~~(i=1,\cdots,12),\quad$$
For arbitrary $\xi\in\mathbb{R}^{4\times11}\setminus\{0\},$ we
have
\begin{equation*}\sum_{ijkl}a_{ikjl}\xi_k^i\xi_l^j=\sum_{kl}\sum_{i}g^{kl}\xi_k^i\xi_l^i>0.
\end{equation*}

Summarize  the above discussions, we have the following lemma.

\begin{Lemma}  The
differential operator of the right hand side of (\ref{20}) with
respect to the metric $g$ is uniformly elliptic.\end{Lemma}

\noindent{\bf Proof of Theorem 4.1.} Noting Lemmas 4.2, 4.4, 4.5
and the compactness property of $M,$ and using the standard
theorem of partial differential equations (see  \cite{Amann1},
\cite{Amann2}, \cite{Giaquinta}), we can immediately obtain the
local existence of smooth solution of the modified system
(\ref{20}) with the initial value
$$g_{ij}(x,0)=\tilde g_{ij}(x),\quad A_i(x,0)=\tilde A_i(x),\quad
B_{ij}(x,0)=\tilde B_{ij}(x).$$ In turn the solution of the
gradient flow (\ref{21}) can be obtained from (\ref{19}) (or
(\ref{17})). The proof of the existence of smooth solution is
completed.

Now we argue the uniqueness of the solution of the gradient flow
(\ref{21}).

By Lemma 4.2, 4.4 and the standard theorem of partial differential
equations, we can obtain the uniqueness of $A$ and $B$ . For any
two solutions $\hat g^{(1)}_{ij}$ and $\hat g^{(2)}_{ij} $ of the
gradient flow (\ref{21}) with the same initial data, we can solve
the initial value problem (\ref{19}) (or (\ref{17})) to get two
families $\varphi^{(1)}$ and $\varphi^{(2)}$ of diffeomorphisms of
$M$. Thus we get two solutions
$$g^{(1)}_{ij}(\cdot,t)=(\varphi^{(1)}_t)^*\hat
g^{(1)}_{ij}(\cdot,t),\qquad
g^{(2)}_{ij}(\cdot,t)=(\varphi^{(2)}_t)^*\hat
g^{(2)}_{ij}(\cdot,t),$$ to the modified evolution (\ref{20})
equations with the same initial value $g_{ij}(x,0)=\tilde
g_{ij}(x)$. The uniqueness result for the strictly parabolic
equation implies that $g^{(1)}_{ij}=g^{(2)}_{ij}.$ Since the
initial value problem (\ref{19}) is clearly a strictly parabolic
system, the corresponding solutions $\varphi^{(1)}$ and
$\varphi^{(2)}$ of (\ref{19}) must agree. Consequently, the
metrics $\hat g^{(1)}_{ij}$ and $\hat g^{(2)}_{ij}$ must agree
also. Thus, we have proved Theorem.$\qquad\Box$

\begin{Remark} we are also able to consider a flow for a
similar functional for a four-dimension manifold
$$S_1 = \int_{M} d^{4}x \sqrt{g} e^{- f} ( \chi + R + 4|\nabla\phi
|^{2}) - \frac{\epsilon_{H}}{2}e^{- f} H \wedge * H -
\epsilon_{F}e^{- f} F \wedge
* F + \frac{e}{2} F \wedge F$$where $e$ is the Euler number $e(\eta)$
 of the bundle $\eta.$ The corresponding flow is given by
$$\begin{cases} \dfrac{\partial g_{ij}}{\partial
t}=-2[R_{ij}+\nabla_i\nabla_j f-\dfrac{1}{4}H_{ikl}H_j^{~kl}-F_i^{~k}F_{jk}],\vspace{2mm}\\
\dfrac{\partial B_{ij}}{\partial
t}=e^{f}\nabla_k(e^{-f}H^k_{~ij}),\vspace{2mm}\\
\dfrac{\partial A_i}{\partial
t}=-e^{f}\nabla_k(e^{-f}F^{~k}_{i}),\vspace{2mm}\\
\dfrac{\partial f}{\partial t}=\chi-2R-3\triangle f+|\nabla
f|^2+\dfrac{1}{3}H^2+\dfrac{3}{2}F^2.
\end{cases}$$By the same argument, we can obtain the same results in
section 3-4.\end{Remark}

\section{The Monotonicity Formula}

Let $M$ be a $n$-dimensional compact Riemannian manifold with
metric $g_{ij}$ , the Levi-Civita connection is given by the
Christoffel symbols
\begin{equation*}
\Gamma^{k}_{ij}=\frac{1}{2}g^{kl}\left\{\frac{\partial
g_{jl}}{\partial x^{i}}+\frac{\partial g_{il}}{\partial
x^{j}}-\frac{\partial g_{ij}}{\partial x^{l}}\right\},
\end{equation*}
where $(g^{ij})$ is the inverse of $(g_{ij})$. The Riemannian
curvature tensors read
\begin{eqnarray*}
R^{k}_{ijl}=\frac{\partial\Gamma^{k}_{jl}}{\partial
x^{i}}-\frac{\partial\Gamma^{k}_{il}}{\partial
x^{j}}+\Gamma^{k}_{ip}\Gamma^{p}_{jl}-\Gamma^{k}_{jp}\Gamma^{p}_{il},\;\;\;
R_{ijkl}=g_{kp}R^{p}_{ijl}.
\end{eqnarray*}
The Ricci tensor is the contraction
\begin{equation*}
R_{ik}=g^{jl}R_{ijkl}
\end{equation*}
and the scalar curvature is
\begin{equation*}
R=g^{ij}R_{ij}.
\end{equation*}

 For each field we shall consider the gauge equivalent classes
of fields. Two metrics $g_{1}, g_{2}$ are in the same equivalent
class if and only if they are differ by a diffeomorphism, i.e.,
there exists a diffeomorphism $f: M \rightarrow M$ such that
$g_{2} = f^{*} g_{1}.$ Two gauge fields $A_{1}$ and $A_{2}$ are
equivalent if and only if there exists a function $\alpha$ on $M$
such that $A_{2} = A_{1} + d \alpha.$ Two $B$-fields $B_{1}$ and
$B_{2}$ are equivalent if and only if there exists an one-form
$\beta$ on $M$ such that $B_{2} = B_{1} + d \beta.$


From the first variation of $S$, we can obtain the flow equations
$$\begin{cases} \dfrac{\partial g_{ij}}{\partial
t}=-2[R_{ij}+\nabla_i\nabla_jf-\dfrac{1}{4}H_{ikl}H_j^{~kl}-F_i^{~k}F_{jk}],\vspace{2mm}\\
\dfrac{\partial B_{ij}}{\partial
t}=e^{f}\nabla_k(e^{-f}H^k_{~ij}),\vspace{2mm}\\
\dfrac{\partial A_i}{\partial
t}=-e^{f}\nabla_k(e^{-f}F^{~k}_{i}),\vspace{2mm}\\
\dfrac{\partial \phi}{\partial t}=\chi-2R-3\triangle f+|\nabla
f|^2+\dfrac{1}{3}H^2+\dfrac{3}{2}F^2.
\end{cases}$$
If $\varphi_t$ is a one-parameter group of diffeomorphisms
generated by a vector field $\nabla f$, we have
\begin{equation*}\label{14}\dfrac{\partial g_{ij}}{\partial
t}=-2(R_{ij}-\dfrac14H_{ikl}H_{j}^{~kl}-F_{ik}F_{j}^{~k}),\end{equation*}
\begin{equation*}\dfrac{\partial A_i}{\partial t}=-
\nabla_kF_{i}^{~k}+\dfrac{\partial}{\partial
x^i}(\nabla^kfA_k),\end{equation*}
\begin{equation*}\dfrac{\partial B_{ij}}{\partial t}=
\nabla_kH_{~ij}^{k}+\dfrac{\partial}{\partial
x^i}(\nabla^kfB_{kj})+\dfrac{\partial}{\partial
x^j}(\nabla^kfB_{ik}).\end{equation*} Let $\tilde{A}=A-d\beta$
where $\dfrac{\partial\beta}{\partial t }=\nabla^kfA_k,$ then
$\tilde{F}=F$ and
\begin{equation*}\dfrac{\partial\tilde {A}_i}{\partial t}=-\nabla_k\tilde {F}_{i}^{~k}.\end{equation*}
Similarly, let $\tilde{B}=B+d\omega$ where
$\dfrac{\partial\omega_i}{\partial t}=\nabla^kfB_{ik},$ then
\begin{equation*}\dfrac{\partial\tilde{B}_{ij}}{\partial t}=\nabla_k(\tilde{H}^{k}_{~ij}).\end{equation*}

Because $A$ and $\tilde A$ ($B$ and $\tilde B$) are in the same
gauge equivalent class, we still denote $\tilde A$ ($\tilde B$) as
$A$ ($B$). Now we consider the flow equation
\begin{equation}\label{14}\dfrac{\partial g_{ij}}{\partial
t}=-2(R_{ij}-\dfrac14H_{ikl}H_{j}^{~kl}-F_{ik}F_{j}^{~k}),\end{equation}
\begin{equation}\label{15}\dfrac{\partial A_i}{\partial t}=-\nabla_k {F}_{i}^{~k},\end{equation}
\begin{equation}\label{16}\dfrac{\partial{B}_{ij}}{\partial t}=\nabla_k({H}^{k}_{~ij}).\end{equation}

\begin{Theorem}Let $g_{ij}, A_i, B_{ij} ~\text{and}~ f$ evolve according to the coupled flow
$$\begin{cases} \dfrac{\partial g_{ij}}{\partial
t}=-2[R_{ij}-\dfrac14H_{ikl}H_j^{~kl}-F_i^{~k}F_{jk}],\vspace{2mm}\\
\dfrac{\partial B_{ij}}{\partial
t}=\nabla_k{H}^{k}_{~ij},\vspace{2mm}\\
\dfrac{\partial A_i}{\partial
t}=-\nabla_k {F}_{i}^{~k},\vspace{2mm}\\
\dfrac{\partial f}{\partial t}=\chi-2R-3\triangle f+2|\nabla
f|^2+\dfrac13H^2+\dfrac32F^2.\end{cases}$$ Then
$$\begin{aligned}\dfrac{d S}{dt}&=\int\left[(-\chi+R-|\nabla f|^2+2\triangle f-\dfrac{1}{12}H^2-
\dfrac12F^2)^2+2(R_{ij}+\nabla_i\nabla_j f
-\dfrac14H_{ikl}H_j^{~kl}-F_i^{~k}F_{jk})^2\right.\\&\left.~~~+2(\nabla_kF_{i}^{~k}-F_{i}^{~k}\nabla_kf)^2+\dfrac12(\nabla_kH_{k}
^{~ij}-H_{k} ^{~ij}\nabla_kf)^2\right]e^{-f}d V.\end{aligned}$$
 In particular $S$ is nondecreasing in time and the
monotonicity is strict unless we are on the critical points.
\end{Theorem}

\noindent{\bf Proof.}
\begin{equation*}\label{2}\begin{aligned}\dfrac{d S}{dt}&=\int d^3x\sqrt{g}e^{-f}(\dfrac12g^{ij}\dfrac
 {\partial g_{ij}}{\partial t}-\dfrac
 {\partial f}{\partial t})
(-\chi+R+2\triangle f-|\nabla
f|^2-\dfrac{1}{12}H^2-\dfrac{1}{2}F^2)\\
&~~~+\int d^3x\sqrt{g}e^{-f}\dfrac
 {\partial g_{ij}}{\partial t}(-R_{ij}-\nabla_i\nabla_j
f+\dfrac14H_{ikl}H_j^{~kl}+F_i^{~k}F_{jk})\\
&~~~+\int d^3x\sqrt{g}e^{-f}\dfrac
 {\partial A_i}{\partial t}(-2\nabla_k(F_i^{~k}e^{-f})e^f)+\dfrac
 {\partial B_{ij}}{\partial t}(\dfrac12\nabla_k(H^k_{~ij}e^{-f})e^f)\\
&=\int(\triangle f-|\nabla f|^2)(-\chi+R-|\nabla f|^2+2\triangle
f-\dfrac{1}{12}H^2-\dfrac12F^2)e^{-f}d
V\\&~~~+\int[-\chi+R-|\nabla
f|^2+2\triangle f-\dfrac{1}{12}H^2-\dfrac12F^2]^2e^{-f}d V\\
&~~~+\int2(R_{ij}+\nabla_i\nabla_jf-\dfrac14H_{ikl}H^{~kl}_{j}-F_{ik}F^{~k}_{j})^2e^{-f}d
V\\&~~~+\int-2\nabla_i\nabla_jf(R_{ij}+\nabla_i\nabla_jf-
\dfrac14H_{ikl}H^{~kl}_{j}-F_{ik}F^{~k}_{j})e^{-f}d V\\
&~~~+\int2(\nabla_kF_{i}^{~k}-F_{i}^{~k}\nabla_kf)^2e^{-f}d
V+\int\dfrac12(\nabla_kH^{k}_{~ij}-H^{k}_{~ij}\nabla_kf)^2e^{-f}d V\\
&~~~+\int2F_{i}^{~k}\nabla_kf(\nabla_kF_{i}^{~k}-F_{i}^{~k}\nabla_kf)e^{-f}d
V+\int\dfrac12H^{k}_{~ij}\nabla_kf(\nabla_kH^{k}_{~ij}-H^{k}_{~ij}\nabla_kf)e^{-f}d
V .\end{aligned}\end{equation*} By the similar argument of Ricci
flow, we have
$$\int(\triangle f-|\nabla f|^2)(R-|\nabla f|^2+2\triangle f)e^{-f}d V=2\int\nabla_i\nabla_jf(\nabla_i\nabla_jf+R_{ij}
)e^{-f}d V.$$ And noting the following properties
$$\nabla_mF_{ij}+\nabla_jF_{mi}+\nabla_iF_{jm}=0,$$
$$\nabla_mH_{ijk}=\nabla_iH_{mjk}+\nabla_jH_{imk}+\nabla_kH_{ijm},$$
we have
\begin{equation*}\begin{aligned}&~~~\int(\triangle f-|\nabla f|^2)(-\chi-\dfrac{1}{12}H^2-\dfrac12F^2)
e^{-f}d V\\
&=\int
g^{ij}(\nabla_i\nabla_jf-\nabla_if\nabla_jf)(-\chi-\dfrac{1}{12}H^2-\dfrac12F^2)
e^{-f}d V\\&=\int g^{ij}\nabla_if
\nabla_j(\chi+\dfrac{1}{12}H^2+\dfrac12F^2)e^{-f}d V\\
&=\int g^{ij}\nabla_if(\dfrac16\nabla_jH_{pkl}H^{pkl}+\nabla_jF_{kl}F^{kl})e^{-f}dV\\
&=\int
g^{ij}\nabla_if(\dfrac16(\nabla_pH_{jkl}+\nabla_kH_{pjl}+\nabla_lH_{pkj})H^{pkl}+(-\nabla_kF_{lj}-
\nabla_lF_{jk})F^{kl})e^{-f}dV\\
&=\int g^{ij}\nabla_if(\dfrac12\nabla_pH_{jkl}H^{pkl}+2\nabla_kF_{jl}F^{kl})e^{-f}dV\\
&=\int(-\dfrac12g^{ij}\nabla_p\nabla_ifH_{jkl}H^{pkl}-2g^{ij}\nabla_k\nabla_ifF_{jl}F^{kl})e^{-f}dV\\
&~~~+\int\dfrac12g^{ij}H_{jkl}\nabla_if(-\nabla_pH^{pkl}+\nabla_pfH^{pkl})e^{-f}dV
+\int2g^{ij}\nabla_ifF_{jl}(\nabla_kfF^{kl}-\nabla_kF^{kl})e^{-f}dV\\
&=\int
2\nabla_i\nabla_jf(-\dfrac14H_{ikl}H_{j}^{~kl}-F_{ik}F^{~k}_{j})e^{-f}dV+\int\dfrac12\nabla_kfH_{~ij}^{k}
(H^{p}_{~ij}\nabla_pf-\nabla_pH_{~ij}^{p})e^{-f}dV\\
&~~~+\int2\nabla_kfF_i^{~k}(F_i^{~k}\nabla_kf-\nabla_kF_{i}^{~k})e^{-f}dV.\end{aligned}\end{equation*}
Combining with the above argument, we finish the proof.

Let $u=e^{-f}$ be the lowest eigenfunction of the Schrodinger
operator, i.e.

$$(R - \frac{1}{12} H^{2} - \frac{1}{2} F^{2} - 4 \Delta) u = \lambda u,$$

or,

$$R - \frac{1}{12} H^{2} - \frac{1}{2} F^{2} + 2 \Delta f - |\nabla f|^{2} =
\lambda.$$

It minimizes the functional

$$S(g, A, B, f) = \int_{M} dV e^{-f/2} (R - \frac{1}{12} H^{2} - \frac{1}{2} F^{2}
- 4 \Delta) e^{- f/2} / \int_{M} e^{-f} dV.$$

We have a new functional

$$\lambda(g, A, B) = inf_{\{ f | \int_{M} e^{-f} dV = 1\}} S(g, A,
B, f).$$

Let $\lambda(t) = \lambda(g(t), A(t), B(t))$, we have

$$\frac{d \lambda}{dt} = \int_{M} (|R_{ij} + \nabla_{i} \nabla_{j} f -
\frac{1}{4} H_{ikj} H^{kl}_{j} - F_{ik} F^{k}_{j}|^{2} +
\frac{1}{4} |\nabla^{k} H_{kij} - H_{kij} \nabla^{k}f|^{2} +
|\nabla_{k} F^{k}_{i} - F^{k}_{i} \nabla_{k} f |^{2} ) e^{-f}
dV.$$

We have then (see also \cite{Oliynyk}):

1) $\lambda(t)$ is monotone, i.e. $\frac{d \lambda(t)}{dt} \ge 0$.

2) Critical points of (*) are the same as critical points of
$\lambda$.

\section{Critical points}

Consider the functional
\begin{equation}\label{1}\begin{aligned}S &= \int_{M}
d^{3}x \sqrt{g} e^{- f} ( -\chi + R + |\nabla f|^{2}) -
\frac{1}{2}e^{- f} H \wedge
* H -
e^{- f} F \wedge * F \\
&=\int d^3x\sqrt{g}e^{-f}(-\chi+R+|\nabla
f|^2-\dfrac{1}{12}H^2-\dfrac12F^2).\end{aligned}\end{equation} Its
first variation can be expressed as follows
\begin{equation}\label{2}\begin{aligned}\delta S=&\int d^3x\sqrt{g}e^{-f}(\dfrac12g^{ij}\delta g_{ij}-\delta f)
(-\chi+R+2\triangle f-|\nabla
f|^2-\dfrac{1}{12}H^2-\dfrac{1}{2}F^2)\\
&+\int d^3x\sqrt{g}e^{-f}\delta g_{ij}(-R_{ij}-\nabla_i\nabla_j
f+\dfrac14H_{ikl}H_j^{~kl}+F_i^{~k}F_{jk})\\
&+\int d^3x\sqrt{g}e^{-f}\delta
A_i(-2\nabla_k(F_i^{~k}e^{-f})e^f)+\delta
B_{ij}(\dfrac12\nabla_k(H^k_{~ij}e^{-f})e^f).\end{aligned}\end{equation}
The $U(1)$ gauge field $A$ is a one-form potential whose field
strength $F=dA.$ The Wess-Zumino field $B$ is a two-form potential
whose field strength $H = dB , \eta$ is the volume form, $f$ is a
dilaton. And in 3-dimension manifold, the field strength is
proportional to the Levi-Civita tensor
$H_{\mu\nu\rho}=H(x)\eta_{\mu\nu\rho},$ where $H(x)$ is a scalar
field and $\eta^{\mu\nu\rho}=\epsilon^{\mu\nu\rho}/\sqrt g$ is the
completely skewsymmetric  Levi-Civita tensor. Therefore, the
critical points satisfy the following equations
\begin{equation}\label{3}R_{ij}+\nabla_i\nabla_j
f-\dfrac14H_{ikl}H_j^{~kl}-F_i^{~k}F_{jk}=0,\end{equation}
\begin{equation}\label{4}\nabla_k(F_i^{~k}e^{-f})=0,\end{equation}
\begin{equation}\label{5}\nabla_k(H^k_{~ij}e^{-f})=0,\end{equation}
\begin{equation}\label{6}-\chi+R+2\triangle f-|\nabla
f|^2-\dfrac{1}{12}H^2-\dfrac{1}{2}F^2=0.\end{equation}

Suppose $M$ is a compact Riemannian manifold. From (\ref{4}) and
(\ref{5}), we can obtain $F=H=0$ at the critical points of the
general Ricci flow on $M$. In fact,
\begin{equation*}\begin{aligned}\int_MF^2e^{-f}dV&=\int_MF^{ij}F_{ij}e^{-f}dV=\int_MF^{ij}
(\nabla_iA_j-\nabla_jA_i)e^{-f}dV\\ &=
2\int_MF^{ij}\nabla_iA_je^{-f}dV=-2\int_M\nabla_i(F^{ij}e^{-f})A_jdV=0,\end{aligned}\end{equation*}

\begin{equation*}\begin{aligned}\int_MH^2e^{-f}dV&=\int_MH^{ijk}H_{ijk}e^{-f}dV=\int_MH^{ijk}
(\nabla_kB_{ij}+\nabla_iB_{jk}+\nabla_jB_{ki})e^{-f}dV\\ &=
3\int_MH^{ijk}\nabla_iB_{jk}e^{-f}dV=-3\int_M\nabla_i(H^{ijk}e^{-f})B_{jk}dV=0.\end{aligned}\end{equation*}

{\bf Remark:} Although the fields $F$ and $H$ do not provide any
help in the study of critical points of general Ricci flow for
compact Riemannian manifold, they maybe play an important role for
the noncompact case.

\section{Evolution of Curvatures}

By virtue of the curvature tensor evolution equations of the Ricci
flow, we can obtain the curvature tensor evolution equations under
the gradient flow (\ref{21}). Let us choose a normal coordinate
system $\{x^i\}$ around a fixed point $x\in M$ such that
$\dfrac{\partial g_{ij}}{\partial x^k}=0$ and
$g_{ij}(p)=\delta_{ij}.$

\begin{Theorem}  Under the gradient flow (\ref{21}), the
curvature tensor satisfies the evolution equation
\begin{equation*}\begin{aligned}
\dfrac{\partial}{\partial t}R_{ijkl}&=\triangle
R_{ijkl}+2(B_{ijkl}-B_{ijlk}-B_{iljk}+B_{ikjl})-g^{pq}(R_{pjkl}R_{qi}+R_{ipkl}R_{qj}+R_{ijpl}R_{qk}+R_{ijkp}R_{ql})\\
&~~~+\dfrac{1}{4}[\nabla_i\nabla_l(H_{kpq}H_j^{~pq})-\nabla_i\nabla_k(H_{jpq}H_l^{~pq})-
\nabla_j\nabla_l(H_{kpq}H_i^{~pq})+\nabla_j\nabla_k(H_{ipq}H_l^{~pq})]\\
&~~~+\dfrac{1}{4}g^{mn}(H_{kpq}H_m^{~pq}R_{ijnl}+H_{mpq}H_l^{~pq}R_{ijkn})\\
&~~~+\nabla_i\nabla_l(F_k^{~p}F_{jp})-\nabla_i\nabla_k(F_j^{~p}F_{lp})-
\nabla_j\nabla_l(F_k^{~p}F_{ip})+\nabla_j\nabla_k(F_i^{~p}F_{lp})\\
&~~~+g^{mn}(F_k^{~p}F_{mp}R_{ijnl}+
F_m^{~p}F_{lp}R_{ijkn}),\\
\end{aligned}\end{equation*}
where $B_{ijkl}=g^{pr}g^{qs}R_{piqj}R_{rksl}$ and $\triangle$ is
the Laplacian with respect to the evolving metric.\end{Theorem}

\noindent{\bf Proof.} At the point $x\in M$, which we has chosen a
normal coordinate system such that $\dfrac{\partial
g_{ij}}{\partial x^k}=0,$  we compute
$$\dfrac{\partial}{\partial t}\Gamma_{jl}^h=\frac12\dfrac{\partial}{\partial t}g^{hm}\left(\dfrac{\partial g_{ml}}
{\partial x^j}+\dfrac{\partial g_{mj}} {\partial
x^l}-\dfrac{\partial g_{jl}} {\partial
x^m}\right)+\frac12g^{hm}\left[\dfrac{\partial} {\partial
x^j}(\dfrac{\partial g_{ml}} {\partial t})+\dfrac{\partial}
{\partial x^l}(\dfrac{\partial g_{mj}} {\partial
t})-\dfrac{\partial} {\partial x^m}(\dfrac{\partial g_{jl}}
{\partial t})\right],$$
\begin{equation*}\begin{aligned}\dfrac{\partial}{\partial t}R_{ijl}^h&=\dfrac{\partial}{\partial x^i}
\left(\dfrac{\partial}{\partial
t}\Gamma_{jl}^h\right)-\dfrac{\partial}{\partial
x^j}\left( \dfrac{\partial}{\partial t}\Gamma_{il}^h\right)\vspace{2mm}\\
&=-\frac12g^{hp}g^{qm}\dfrac{\partial g_{pq}}{\partial t}
\left(\dfrac{\partial^2 g_{ml}}{\partial x^i\partial
x^j}+\dfrac{\partial^2 g_{mj}}{\partial x^i\partial
x^l}-\dfrac{\partial^2 g_{jl}}{\partial x^i\partial
x^m}\right)\vspace{2mm}\\
&~~~+\frac12g^{hp}g^{qm}\dfrac{\partial g_{pq}}{\partial t}
\left(\dfrac{\partial^2 g_{ml}}{\partial x^j\partial
x^i}+\dfrac{\partial^2 g_{mi}}{\partial x^j\partial
x^l}-\dfrac{\partial^2 g_{il}}{\partial x^j\partial
x^m}\right)\vspace{2mm}\\
&~~~+\frac12g^{hm}\left[\dfrac{\partial^2 }{\partial x^i\partial
x^l}\left(\dfrac{\partial g_{mj}} {\partial
t}\right)-\dfrac{\partial^2 }{\partial x^i\partial
x^m}\left(\dfrac{\partial g_{jl}} {\partial
t}\right)-\dfrac{\partial^2 }{\partial x^j\partial
x^l}\left(\dfrac{\partial g_{mi}} {\partial
t}\right)+\dfrac{\partial^2 }{\partial x^j\partial
x^m}\left(\dfrac{\partial g_{il}} {\partial t}\right)\right]\vspace{2mm}\\
&=\frac12g^{hm}\left[\dfrac{\partial^2 }{\partial x^i\partial
x^l}\left(\dfrac{\partial g_{mj}} {\partial
t}\right)-\dfrac{\partial^2 }{\partial x^i\partial
x^m}\left(\dfrac{\partial g_{jl}} {\partial
t}\right)-\dfrac{\partial^2 }{\partial x^j\partial
x^l}\left(\dfrac{\partial g_{mi}} {\partial
t}\right)+\dfrac{\partial^2 }{\partial x^j\partial
x^m}\left(\dfrac{\partial g_{il}} {\partial t}\right)\right]\vspace{2mm}\\
&~~~-g^{hp}\dfrac{\partial g_{pq}}{\partial t}R_{ijl}^q\,\,,\vspace{2mm}\\
\dfrac{\partial}{\partial t}R_{ijkl}& =\dfrac{\partial}{\partial
t}R_{ijl}^hg_{kh}+R_{ijl}^h\dfrac{\partial}{\partial t}g_{kh}\\
&=\frac12\left[\dfrac{\partial^2 }{\partial x^i\partial
x^l}\left(\dfrac{\partial g_{kj}} {\partial
t}\right)-\dfrac{\partial^2 }{\partial x^i\partial
x^k}\left(\dfrac{\partial g_{jl}} {\partial
t}\right)-\dfrac{\partial^2 }{\partial x^j\partial
x^l}\left(\dfrac{\partial g_{ki}} {\partial
t}\right)+\dfrac{\partial^2 }{\partial x^j\partial
x^k}\left(\dfrac{\partial g_{il}} {\partial
t}\right)\right]\ \ ,\\
\end{aligned}\end{equation*}
then we have
\begin{equation*}\begin{aligned}\dfrac{\partial}{\partial t}R_{ijkl}
&=\dfrac{\partial^2 }{\partial x^i\partial
x^k}R_{jl}-\dfrac{\partial^2 }{\partial x^i\partial
x^l}R_{kj}+\dfrac{\partial^2 }{\partial x^j\partial
x^l}R_{ki}-\dfrac{\partial^2 }{\partial x^j\partial x^k}R_{il}\\
&~~~+\dfrac{1}{4}\left[\dfrac{\partial^2 }{\partial x^i\partial
x^l}(H_{kpq}H_j^{~pq})-\dfrac{\partial^2 }{\partial x^i\partial
x^k}(H_{jpq}H_l^{~pq})-\dfrac{\partial^2 }{\partial x^j\partial
x^l}(H_{kpq}H_i^{~pq})+\dfrac{\partial^2 }{\partial x^j\partial
x^k}(H_{ipq}H_l^{~pq})\right]\\
&~~~+\dfrac{\partial^2 }{\partial x^i\partial
x^l}(F_k^{~p}F_{jp})-\dfrac{\partial^2 }{\partial x^i\partial
x^k}(F_j^{~p}F_{lp})-\dfrac{\partial^2 }{\partial x^j\partial
x^l}(F_k^{~p}F_{ip})+\dfrac{\partial^2 }{\partial x^j\partial
x^k}(F_i^{~p}F_{lp})\\
&\triangleq I_1+\dfrac{1}{4}I_2+I_3.
\end{aligned}\end{equation*}

By the identity (see \cite{Cao})
\begin{equation*}\begin{aligned}&~~~\nabla_i\nabla_k
R_{jl}-\nabla_i\nabla_l R_{jk}-\nabla_j\nabla_k
R_{il}+\nabla_j\nabla_l R_{ik}\\
&=\triangle
R_{ijkl}+2(B_{ijkl}-B_{ijlk}-B_{iljk}+B_{ikjl})-g^{pq}(R_{pjkl}R_{qi}+R_{ipkl}R_{qj})
\end{aligned}\end{equation*}and
$$\nabla_i\nabla_k
R_{jl}=\dfrac{\partial^2 R_{jl}}{\partial x^i\partial
x^k}-R_{ml}\dfrac{\partial}{\partial
x^i}\Gamma^m_{kj}-R_{jm}\dfrac{\partial}{\partial
x^i}\Gamma^m_{kl}\,\,,$$ we have
\begin{equation*}\begin{aligned}I_1&=\nabla_i\nabla_k
R_{jl}+R_{ml}\dfrac{\partial}{\partial
x^i}\Gamma^m_{kj}+R_{jm}\dfrac{\partial}{\partial
x^i}\Gamma^m_{kl}-\nabla_i\nabla_l
R_{kj}-R_{km}\dfrac{\partial}{\partial
x^i}\Gamma^m_{lj}-R_{mj}\dfrac{\partial}{\partial
x^i}\Gamma^m_{lk}\\
&~~~-\nabla_j\nabla_k R_{il}-R_{ml}\dfrac{\partial}{\partial
x^j}\Gamma^m_{ki}-R_{im}\dfrac{\partial}{\partial
x^j}\Gamma^m_{kl}+\nabla_j\nabla_l
R_{ki}+R_{km}\dfrac{\partial}{\partial
x^j}\Gamma^m_{li}+R_{mi}\dfrac{\partial}{\partial
x^j}\Gamma^m_{lk}\\
&=\nabla_i\nabla_k R_{jl}-\nabla_i\nabla_l R_{jk}-\nabla_j\nabla_k
R_{il}+\nabla_j\nabla_l
R_{ik}-R_{km}R^m_{ijl}+R_{ml}R_{ijk}^m\\
&=\triangle
R_{ijkl}+2(B_{ijkl}-B_{ijlk}-B_{iljk}+B_{ikjl})-g^{pq}(R_{pjkl}R_{qi}+R_{ipkl}R_{qj}+R_{ijpl}R_{qk}+R_{ijkp}R_{ql}),\\
\end{aligned}\end{equation*}where $B_{ijkl}=g^{pr}g^{qs}R_{piqj}R_{rksl}.$

Now we compute $I_2.$

It is easily verified that
$$\nabla_i\nabla_k(H_{jpq}H^{~pq}_{l})=\dfrac{\partial^2}{\partial x^i\partial x^k}(H_{jpq}H^{~pq}_{l})-
H_{mpq}H^{~pq}_{l}\dfrac{\partial}{\partial
x^i}\Gamma^m_{kj}-H_{jpq}H^{~pq}_{m}\dfrac{\partial}{\partial
x^i}\Gamma^m_{kl}\,\,.$$ As a result, we obtain
\begin{equation*}\begin{aligned}I_2&=\nabla_i\nabla_l(H_{kpq}H^{~pq}_{j})+
H_{mpq}H^{~pq}_j\dfrac{\partial}{\partial
x^i}\Gamma^m_{lk}+H_{kpq}H^{~pq}_{m}\dfrac{\partial}{\partial
x^i}\Gamma^m_{lj}-\nabla_i\nabla_k(H_{jpq}H^{~pq}_{l})\\
&~~~- H_{mpq}H^{~pq}_{l}\dfrac{\partial}{\partial
x^i}\Gamma^m_{kj}-H_{jpq}H^{~pq}_{m}\dfrac{\partial}{\partial
x^i}\Gamma^m_{kl}-\nabla_j\nabla_l(H_{kpq}H^{~pq}_{i})-
H_{mpq}H^{~pq}_{i}\dfrac{\partial}{\partial
x^j}\Gamma^m_{lk}\\
&~~~-H_{kpq}H^{~pq}_{m}\dfrac{\partial}{\partial
x^j}\Gamma^m_{li}+\nabla_j\nabla_k(H_{ipq}H^{~pq}_{l})+
H_{mpq}H^{~pq}_{l}\dfrac{\partial}{\partial
x^j}\Gamma^m_{ki}+H_{ipq}H^{~pq}_{m}\dfrac{\partial}{\partial
x^j}\Gamma^m_{kl}\\
&=\nabla_i\nabla_l(H_{kpq}H^{~pq}_{j})-\nabla_i\nabla_k(H_{jpq}H^{~pq}_{l})-\nabla_j\nabla_l(H_{kpq}H^{~pq}_{i})
+\nabla_j\nabla_k(H_{ipq}H^{~pq}_{l})\\
&~~~+H_{kpq}H^{~pq}_{m}R^m_{ijl}+H_{mpq}H^{~pq}_{l}R^m_{jik}\\
&=\nabla_i\nabla_l(H_{kpq}H^{~pq}_{j})-\nabla_i\nabla_k(H_{jpq}H^{~pq}_{l})-\nabla_j\nabla_l(H_{kpq}H^{~pq}_{i})
+\nabla_j\nabla_k(H_{ipq}H^{~pq}_{l})\\
&~~~+g^{mn}(H_{kpq}H^{~pq}_{m}R_{ijnl}+H_{mpq}H^{~pq}_{l}R_{ijkn})\,\,.\\
\end{aligned}\end{equation*}
Now it remains to compute the last term.  The following identity
$$\nabla_i\nabla_k(F_j^{~p}F_{lp})=\dfrac{\partial^2}{\partial x^i\partial x^k}(F_j^{~p}F_{lp})-
F_m^{~p}F_{lp}\dfrac{\partial}{\partial
x^i}\Gamma^m_{kj}-F_j^{~p}F_{mp}\dfrac{\partial}{\partial
x^i}\Gamma^m_{kl}\,\,$$  yields
\begin{equation*}\begin{aligned}I_3&=\nabla_i\nabla_l(F_k^{~p}F_{jp})+
F_m^{~p}F_{jp}\dfrac{\partial}{\partial
x^i}\Gamma^m_{lk}+F_k^{~p}F_{mp}\dfrac{\partial}{\partial
x^i}\Gamma^m_{lj}-\nabla_i\nabla_k(F_j^{~p}F_{lp})\\
&~~~- F_m^{~p}F_{lp}\dfrac{\partial}{\partial
x^i}\Gamma^m_{kj}-F_j^{~p}F_{mp}\dfrac{\partial}{\partial
x^i}\Gamma^m_{kl}-\nabla_j\nabla_l(F_k^{~p}F_{ip})-
F_m^{~p}F_{ip}\dfrac{\partial}{\partial
x^j}\Gamma^m_{lk}\\
&~~~-F_k^{~p}F_{mp}\dfrac{\partial}{\partial
x^j}\Gamma^m_{li}+\nabla_j\nabla_k(F_i^{~p}F_{lp})+
F_m^{~p}F_{lp}\dfrac{\partial}{\partial
x^j}\Gamma^m_{ki}+F_i^{~p}F_{mp}\dfrac{\partial}{\partial
x^j}\Gamma^m_{kl}\\
&=\nabla_i\nabla_l(F_k^{~p}F_{jp})-\nabla_i\nabla_k(F_j^{~p}F_{lp})-\nabla_j\nabla_l(F_k^{~p}F_{ip})
+\nabla_j\nabla_k(F_i^{~p}F_{lp})\\
&~~~+g^{mn}(F_k^{~p}F_{mp}R_{ijnl}+F_m^{~p}F_{lp}R_{ijkn})\,\,.\\
\end{aligned}\end{equation*}
Combining the above discussions, we complete the proof of the
theorem.$\qquad\Box$

\begin{Theorem} The Ricci curvature satisfies the following evolution equation
\begin{equation*}\begin{aligned}\dfrac{\partial}{\partial
t}R_{ik}&=\triangle
R_{ik}+2g^{pr}g^{qs}R_{piqk}R_{rs}-2g^{pq}R_{pi}R_{qk}\\
&~~~+\dfrac{1}{4}g^{jl}[\nabla_i\nabla_l(H_{kpq}H_j^{~pq})-\nabla_i\nabla_k(H_{jpq}H_l^{~pq})-
\nabla_j\nabla_l(H_{kpq}H_i^{~pq})+\nabla_j\nabla_k(H_{ipq}H_l^{~pq})]\\
&~~~+\dfrac{1}{4}g^{mn}(H_{kpq}H_m^{~pq}R_{in}-g^{jl}H_{mpq}H_l^{~pq}R_{ijkn})\\
&~~~+g^{jl}[\nabla_i\nabla_l(F_k^{~p}F_{jp})-\nabla_i\nabla_k(F_j^{~p}F_{lp})-
\nabla_j\nabla_l(F_k^{~p}F_{ip})+\nabla_j\nabla_k(F_i^{~p}F_{lp})]\\
&~~~+g^{mn}(F_k^{~p}F_{mp}R_{in}-g^{jl}
F_m^{~p}F_{lp}R_{ijkn}).\\
\end{aligned}\end{equation*}
\end{Theorem}
\noindent{\bf Proof.} By Theorem 5.1, we can compute
\begin{equation*}\begin{aligned}\dfrac{\partial}{\partial
t}R_{ik}&=\dfrac{\partial}{\partial
t}R_{ijkl}g^{jl}+R_{ijkl}\dfrac{\partial}{\partial
t}g^{jl}=\dfrac{\partial}{\partial
t}R_{ijkl}g^{jl}-g^{jp}g^{lq}R_{ijkl}\dfrac{\partial}{\partial
t}g_{pq}\\
&=g^{jl}[\triangle
R_{ijkl}+2(B_{ijkl}-B_{ijlk}-B_{iljk}+B_{ikjl})-g^{pq}(R_{pjkl}R_{qi}+R_{ipkl}R_{qj}+R_{ijpl}R_{qk}\\
&~~~+R_{ijkp}R_{ql})]+2g^{jp}g^{lq}R_{ijkl}R_{pq}\\
&~~~+\dfrac{1}{4}g^{jl}[\nabla_i\nabla_l(H_{kpq}H_j^{~pq})-\nabla_i\nabla_k(H_{jpq}H_l^{~pq})-
\nabla_j\nabla_l(H_{kpq}H_i^{~pq})+\nabla_j\nabla_k(H_{ipq}H_l^{~pq})]\\
&~~~+\dfrac{1}{4}g^{jl}g^{mn}(H_{kpq}H_m^{~pq}R_{ijnl}+H_{mpq}H_l^{~pq}R_{ijkn})
-\dfrac{\epsilon_H}{2}R_{ijkl}g^{jp}g^{lq}H_{pmn}H^{~mn}_q\\
&~~~+g^{jl}[\nabla_i\nabla_l(F_k^{~p}F_{jp})-\nabla_i\nabla_k(F_j^{~p}F_{lp})-
\nabla_j\nabla_l(F_k^{~p}F_{ip})+\nabla_j\nabla_k(F_i^{~p}F_{lp})]\\
&~~~+g^{mn}(F_k^{~p}F_{mp}R_{in}+g^{jl}
F_m^{~p}F_{lp}R_{ijkn})-2\epsilon_FR_{ijkl}g^{jp}g^{lq}F_p^{~m}F_{qm}\\
&=\triangle
R_{ik}+2g^{pr}g^{qs}R_{piqk}R_{rs}-2g^{pq}R_{pi}R_{qk}\\
&~~~+\dfrac{1}{4}g^{jl}[\nabla_i\nabla_l(H_{kpq}H_j^{~pq})-\nabla_i\nabla_k(H_{jpq}H_l^{~pq})-
\nabla_j\nabla_l(H_{kpq}H_i^{~pq})+\nabla_j\nabla_k(H_{ipq}H_l^{~pq})]\\
&~~~+\dfrac{1}{4}g^{mn}(H_{kpq}H_m^{~pq}R_{in}-g^{jl}H_{mpq}H_l^{~pq}R_{ijkn})\\
&~~~+g^{jl}[\nabla_i\nabla_l(F_k^{~p}F_{jp})-\nabla_i\nabla_k(F_j^{~p}F_{lp})-
\nabla_j\nabla_l(F_k^{~p}F_{ip})+\nabla_j\nabla_k(F_i^{~p}F_{lp})]\\
&~~~+g^{mn}(F_k^{~p}F_{mp}R_{in}-g^{jl}
F_m^{~p}F_{lp}R_{ijkn}).\qquad\qquad\qquad\qquad\qquad\Box\\
\end{aligned}\end{equation*}

\begin{Theorem}The scalar curvature satisfies the following evolution equation
\begin{equation*}\begin{aligned}\dfrac{\partial}{\partial
t}R&=\triangle
R+2|Ric|^2+\dfrac{1}{2}g^{jl}g^{ik}[\nabla_i\nabla_l(H_{kpq}H^{~pq}_j)
-\nabla_i\nabla_k(H_{jpq}H^{~pq}_l)]\\
&~~~ +2g^{jl}g^{ik}[\nabla_i\nabla_l(F_k^{~p}F_{jp})
-\nabla_i\nabla_k(F_j^{~p}F_{lp})]-g^{ip}R_{ik}\left(\dfrac{1}{2}H_{pmn}H^{kmn}
+2F_{pm}F^{km}\right).\\
\end{aligned}\end{equation*}
\end{Theorem}

\noindent{\bf Proof.} By a direct calculation, we have
\begin{equation*}\begin{aligned}\dfrac{\partial}{\partial
t}R&=\dfrac{\partial}{\partial
t}R_{ik}g^{ik}+R_{ik}\dfrac{\partial}{\partial
t}g^{ik}=\dfrac{\partial}{\partial
t}R_{ik}g^{ik}-R_{ik}g^{ip}g^{kq}\dfrac{\partial}{\partial
t}g_{pq}\\
&=g^{ik}(\triangle
R_{ik}+2g^{pr}g^{qs}R_{piqk}R_{rs}-2g^{pq}R_{pi}R_{qk})+2g^{ip}g^{kq}R_{ik}R_{pq}\\
&~~~+\dfrac{1}{4}g^{jl}g^{ik}[\nabla_i\nabla_l(H_{kpq}H_j^{~pq})-\nabla_i\nabla_k(H_{jpq}H_l^{~pq})-
\nabla_j\nabla_l(H_{kpq}H_i^{~pq})+\nabla_j\nabla_k(H_{ipq}H_l^{~pq})]\\
&~~~+\dfrac{1}{4}g^{ik}g^{mn}(H_{kpq}H_m^{~pq}R_{in}-g^{jl}H_{mpq}H_l^{~pq}R_{ijkn})-\dfrac{1}{2}
g^{ip}g^{kq}R_{ik}H_{pmn}H^{~mn}_q\\
&~~~+g^{ik}g^{jl}[\nabla_i\nabla_l(F_k^{~p}F_{jp})-\nabla_i\nabla_k(F_j^{~p}F_{lp})-
\nabla_j\nabla_l(F_k^{~p}F_{ip})+\nabla_j\nabla_k(F_i^{~p}F_{lp})]\\
&~~~+g^{ik}g^{mn}(F_k^{~p}F_{mp}R_{in}-g^{jl}
F_m^{~p}F_{lp}R_{ijkn})-2g^{ip}g^{kq}R_{ik}F_p^{~m}F_{qm}\\
&=\triangle
R+2|Ric|^2+\dfrac{1}{2}g^{jl}g^{ik}[\nabla_i\nabla_l(H_{kpq}H^{~pq}_j)
-\nabla_i\nabla_k(H_{jpq}H^{~pq}_l)]\\
&~~~ +2g^{jl}g^{ik}[\nabla_i\nabla_l(F_k^{~p}F_{jp})
-\nabla_i\nabla_k(F_j^{~p}F_{lp})]-g^{ip}R_{ik}\left(\dfrac{1}{2}H_{pmn}H^{kmn}
+2F_{pm}F^{km}\right).\  \Box\\
\end{aligned}\end{equation*}

\vskip 5mm

\noindent{\bf Acknowledgements:} The work of S. Hu was supported
in part by the NNSF of China (Grant No. 10771203) and a renovation
grant from the Chinese Academy of Sciences; the work of D. Kong
was supported in part by the NNSF of China (Grant No. 10671124)
and the NCET of China (Grant No. NCET-05-0390); the work of K. Liu
was supported by the NSF and NSF of China.


\begin{thebibliography}{99}

\bibitem{Amann1} H. Amann, {\it Quasilinear evolution equations
and parabolic systems}, Tran. of AMS, {\bf 293} (1986), 191-227.

\bibitem{Amann2} H. Amann, {\it Quasilinear parabolic systems under nonlinear boundary conditions},
Arch. Rat. Mech. Anal., Vol. {\bf 92}  No. {\bf 2} (1986),
153-192.



\bibitem{Cao} H.-D. Cao and X.-P. Zhu, {\it A complete proof of the
Poincar$\acute{e}$ and geometrization conjectures - application of
the Hamilton-Perelman theory of the Ricci flow}, Asian J. Math.,
{\bf 10} (2006), 165-492.

\bibitem{D} D.DeTurck, {\it Deforming metrics in the direction of their Ricci
tensors,} J. Differential Geom., 18 (1983), 157-162.

\bibitem{Gegenberg2} J. Gegenberg and G. Kunstatter, {\it Using 3D stringy gravity to understand the Thurston
conjecture}, arxiv:hep-th/0306279, 2003.

\bibitem{Gegenberg1} J. Gegenberg, S. Vaidya and J. F.
V$\acute{a}$zquez-Poritz, {\it Thurston geometries from eleven
dimensions}, arxiv:hep-th/0205276, 2002.

\bibitem{Giaquinta} M. Giaquinta and G. Modica, {\it Local existence for quasilinear parabolic systems
under nonlinear boundary conditions}, Annali di Matematica Pura ed
Applicata, Vol. {\bf 149} Iss. {\bf 1}(1987), 41-59.

\bibitem{Hamilton} R. S. Hamilton, {\it Three manifolds with positive Ricci
curvature}, J. Differential Geom., {\bf 17} (1982), 255-306.

\bibitem{Oliynyk} T. Oliynyk, V. Suneeta, E. Woolgar, {\it A
gradient flow for worldsheet nonlinear sigma models},
hep-th/0510239.

\bibitem{Perelman} G. Perelman, {\it The entropy formula for the
Ricci flow and its geometric applications}, math.DG/0211159.

\bibitem{Thurston} W. P. Thurston, {\it Three-dimensional geometry and
topology}, Vol.1, Edited by Silvio Levy, Princeton Mathematical
Series, 35, Princeton University Press, Princeton, NJ, 1997.



\end{thebibliography}
\end{document}